\documentstyle[12pt]{amsart}
\begin{document}

\title[{\tiny Preservation of the absolutely continuous spectrum of
Schr\"odinger equation}]
{Preservation of the absolutely continuous spectrum of
Schr\"odinger equation under perturbations by slowly decreasing
potentials and a.e. convergence of integral operators}
\author{Alexander Kiselev}
\address{Mathematical Sciences Research Institute, 1000 Centennial Dr., Berkeley, CA 94720, USA}
\email{kiselev@@msri.org}

\begin{abstract}
We prove new criteria of stability of the absolutely continuous spectrum
of one-dimensional Schr\"odinger operators under slowly decaying
perturbations. As applications, we show that the absolutely continuous 
spectrum of the free and periodic Schr\"odinger operators is preserved 
under perturbations by all potentials $V(x)$ satisfying $|V(x)| \leq 
C(1+x)^{-\frac{2}{3}-\epsilon}.$ The main new technique includes an 
a.e. convergence theorem for a class of integral operators.
\end{abstract}

\maketitle

\begin{center}
 \section*{ Introduction}
\end{center}

\normalsize

In this paper we  study the stability of the absolutely continuous spectrum
of one-dimensional Schr\"odinger operators under perturbations
by slowly decaying potentials. The most general formulation of the
problem we study is as follows. Let $H_{U}$ be a Schr\"odinger
operator acting on $L^{2}(0, \infty)$ and given by the differential expression
\begin{equation}
-\frac{d^{2}}{dx^{2}} + U(x).
\end{equation}
Suppose for simplicity that the potential $U(x)$ is bounded (more general
assumptions on $U(x)$ will be specified below) and fix some self-adjoint
boundary condition at zero.
 It is a well-known fact that under the above
conditions the differential expression (1) defines a unique self-adjoint operator.
Let us assume that the absolutely continuous spectrum of the operator $H_{U}$
is not empty. We perturb the operator $H_{U}$ by a decaying potential $V(x):$
\[ H_{U+V}= H_{U}+V(x). \]
We assume that $V(x)$ is locally integrable and $V(x) \stackrel{x \rightarrow \infty}{\longrightarrow} 0.$
 Since the potential decaying at infinity constitutes a relatively compact
perturbation of the unperturbed part, the essential spectrum of the operator
$H_{U+V}$
coincides with that of $H_{U}.$ We would like to study which conditions on the
rate of decay of $V$ are sufficient to ensure that the absolutely continuous
spectrum of $H_{U}$ is preserved.

The main results of our study is the derivation of the new criteria on the
stability of the absolutely continuous spectrum of Schr\"odinger operators
under perturbations by the general classes of slowly decaying potentials.
We apply these criteria, in particular, to find new classes of slowly decaying
potentials preserving the absolutely continuous spectrum of the free and
periodic Schr\"odinger operators.

Over the years, there has been much attention to the subject and we briefly recall main results.
Suppose that the potential $V$ is of short range, by which we mean that
$|V(x)| \leq C(1+x)^{-1-\epsilon},$ or, more generally, suppose that $V(x)$ is
absolutely integrable. Then for a wide class of background potentials $U(x)$
a number of techniques may be used to show that
the absolutely continuous spectrum of $H_{U}$ is preserved. In many situations the trace class theorems of scattering theory may be used to show that, moreover, the wave operators
exist and are complete in this case. This short-range result is basically all
what is known on the preservation of the absolutely continuous spectrum of
Schr\"odinger operators
under decaying perturbations in general situation.
Essentially more information is available in the case when $U(x)=0.$
 There has been much work on proving the
absolute continuity of the spectrum for the Schr\"odinger operators with
potentials of slower decay satisfying some additional special structural assumptions.
It is a classic result, going back to Weidmann \cite{Weid} that if a
 potential $V$ may be represented as a sum of a function of
bounded variation and an absolutely integrable function, then the spectrum
of the operator $H_{V}$ on $R^{+}=(0, \infty)$ is purely absolutely
contimuous. Many authors developed
a scattering theory for potentials satisfying certain conditions on its
derivatives \cite{Bus}, \cite{As}, \cite{Hor}.  These results hold in any
dimension and the proofs use the method of approximating the scattering trajectories
by the solutions of the classical Hamilton-Jacobi equation. This method is applicable only when there are certian conditions on the decay of the derivatives of the potential. The weakest conditions on the long-range part of the potentials, under which the wave operators exist, are given in \cite{Hor}. For potentials
satisfying $|V(x)| \leq C(1+|x|)^{-\frac{1}{2}-\epsilon},$ for instance, one can
infer the existence of the wave operators if also $|D^{\alpha}V(x)| \leq C_{1}
(1+|x|)^{-\frac{3}{2}-\epsilon}$ for every multiindex $\alpha$ with $|\alpha|=1.$

Another class of results describe spectral behavior of the spherically symmetric (i.e. essentially one-dimensional) specific oscillating potentials, the
typical example being potentials of type $V(x) = \frac{\sin x^{\alpha}}{x^{\beta}}$ with $\alpha,$ $\beta$ positive. Clearly, potentials of this type may not satisfy in general conditions on the decay of the derivatives which are needed for the method of the above works to be applicable. We mention the papers of \cite{Ben},
\cite{BA}, \cite{HiS}, and \cite{White} in which further references may be found. The spectrum of the operator $H_{V}$ for
such potentials turns out to be absolutely continuous with perhaps some isolated
embedded eigenvalues when $\alpha =1.$ Such situations generalize the celebrated
Wigner-von Neumann example \cite{Wig}. Wigner and von Neumann were first 
to construct an example which shows that already for potentials deacying as a Coulomb potential at infinity, i.e. $V(x)=O(\frac{1}{1+|x|}),$
the spectrum need not to be purely absolutely continuous and positive embedded eigenvalues may occur. Naboko \cite{Nab} and  Simon
\cite{Sim1} found different constructions which show that for potentials decaying arbitrarily slower than
at a Coulomb rate, one can already have a rather striking spectral phenomena.
The essential result is that for every function $C(x)$ going monotonically to
infinity as $x$ goes to infinity at an arbitrarily slow rate, and for every
sequence of positive energies $\{\lambda_{i} \}_{i=1}^{\infty}$ one can find
a potential $V(x)$ satisfying
\[ |V(x)| \leq \frac{C(x)}{1+|x|}, \]
such that the corresponding Schr\"odinger operator $H_{V}$ has the sequence
 $\{\lambda_{i}\}^{\infty}_{i=1}$ among its eigenvalues. The earlier work
of Naboko had also an additional assumption of rational independence
of the sequence $\{\sqrt{\lambda_{i}} \}_{i=1}^{\infty},$ but Simon employs
a different method which does away with this condition. In particular, these
examples show that, in general, already for a potential decaying only arbitrary slower than a Coulomb potential, the corresponding Schr\"odinger operator
can have a dense set of eigenvalues on the positive semi-axis. In such situation the Weyl criteria of the preservation of essential spectrum does not tell us
anything about whether there is any other kind of the spectrum on $R^{+}.$
From the constructions of Naboko and Simon, it was also not clear whether the
eigenvalues they construct are really embedded or, in fact, they are the only
spectrum the operator $H_{V}$ has on $R^{+},$ and hence one can have dense
pure point spectrum on $R^{+}$ for the potentials which decay so fast.

 Until recently, even for the basic case $U=0$ the absolutely integrable class
of potentials remained the only class defined purely in terms of the  rate of
decay,  which was known to preserve the absolutely continuous spectrum
of the free Hamiltonian.
A new general class of potentials preserving the absolutely continuous
spectrum of the free Schr\"odinger operator was found in \cite{Kis}.
 Namely, if the potential $V$ satisfies
 $|V(x)| \leq C(1+|x|)^{-\frac{3}{4}-\epsilon}$ with some $\epsilon > 0,$
 with no additional assumptions, then the absolutely continuous spectrum
 fills the whole $R^{+}.$ Of course, as the examples of
Naboko and Simon show, rich embedded singular spectrum may occur; however it
is indeed embedded in the sense that there is an underlying absolutely continuous spectrum. One can decribe the set where the singular part of the spectral
measure might be supported in $R^{+}$ rather explicitly in terms of
the properties of the Fourier transform of $V$ \cite{Kis}.  

We should mention that there exists a work on random potentials by
 Kotani and Ushiroya \cite{Kot} which provides in a sense a bound for the
 best possible result that one
can hope to prove in the deterministic case. They consider the random
potentials of type
\[ V(x) = a(x)F(Y_{x}(\omega)), \]
where $Y_{x}(\omega)$ is a Brownian motion on a compact Riemannian manifold
$M,$ $F$ is a $C^{\infty}$ function which maps $M$ to the real axis,
satisfying certain additional assupmtions (see \cite{Kot}), and $a(x)$ is a
deterministic factor which is in particular taken to be power-decaying:
$a(x)= (1+x)^{-\alpha}.$
The main result of \cite{Kot} shows a sharp transition in the spectral
properties as $\alpha$ passes $\frac{1}{2}:$ for $\alpha < \frac{1}{2},$
the spectrum on the positive half-axis is dense pure point with probability
one; for $\alpha > \frac{1}{2}$ it is almost surely purely absolutely
continuous. For $\alpha =\frac{1}{2}$ one may have a mixture of pure point
and singular continuous spectrum for different regions of energies in
$R^{+}$ with probability one.
Hence, \cite{Kot} implies that there
exist potentials $V(x)$ satisfying  $V(x) \leq C(1+x)^{-\frac{1}{2}},$ which
lead to the purely singular spectrum on $R^{+}.$ Therefore, the largest
general class of power rate decaying potentials for which one can hope to
prove the preservation of the absolutely continuous spectrum of the
free Hamiltonian is the class of potentials satisfying
$|V(x)| \leq C(1+x)^{-\frac{1}{2}-\epsilon}$ with some $\epsilon >0.$

  In this paper we study the problem of the preservation of the absolutely
continuous spectrum of the Schr\"odinger operators in more general setting.
We formulate two general criteria, which may be applied to any operator
$H_{U}$ with non-empty absolutely continuous spectrum, provided it satisfies
a mild additional assumption. These criteria give conditions which are
sufficient for the absolutely continuous spectrum to be stable under perturbations
by all potentials satisfying $|V(x)|\leq C(1+x)^{-\frac{3}{4}-\epsilon}$ and
$|V(x)| \leq C(1+x)^{-\frac{2}{3}-\epsilon}$ respectively.

Shortly, the first criterion works as follows. Suppose that a certain set $S$ belongs to
the essential support of the a.c. part of the spectral measure of the operator $H_{U}.$
Suppose in addition that all solutions of the equation
\begin{equation}
-\frac{d^{2}y}{dx^{2}} + U(x)y = \lambda y.
\end{equation}
 for the energies from
the set $S$ are bounded. This is a rather  natural assumption. In most situations
of interest for almost every energy from the essential support
 of the absolutely continuous spectrum all solutions of (2) are bounded. For every energy $\lambda
 \in S,$ choose two linearly independent uniformly bounded solutions
 $\theta (x, \lambda)$ and $\overline{\theta}(x, \lambda)$ of the equation (2). 
The key object for our study turns out to be the operator $T,$
  defined on the bounded functions of compact support in
 the following way:
\begin{equation} 
 (Tf)(\lambda) = \frac{\chi (S)}{\Im (\theta \overline{\theta}')}
 \int\limits_{0}^{\infty} \theta^{2}(x, \lambda) f(x)\,dx. 
\end{equation}
 Here $\chi (S)$ denotes the characteristic function of the set $S.$
 Note that the factor in front of the integral does not depend on $x$
 since the denominator is the Wronskian of the two linearly
 independent solutions of (2), $\theta$ and
 $\overline{\theta}.$ Suppose that we can pick such
 $\theta (x, \lambda)$ so that the operator $T$
 satisfies an $L^{2}-L^{2}$ bound, i.e. for all measurable bounded functions
 $f$ of compact support we have
 \[ \|Tf\|_{2} \leq C \|f\|_{2}, \]
 where $\|f\|_{2}$ denotes the $L^{2}$ norm of the function $f.$
Then the absolutely continuous spectrum of $H_{U},$ supported on the set
$S,$ is stable under all perturbations by potentials $V(x)$ satisfying
$V(x) \leq C(1+x)^{-\frac{3}{4}-\epsilon}.$ In other words, the set $S$
belongs to the essential support of the absolutely continuous part of the spectral
measure of the Schr\"odinger operator $H_{U+V}.$ We remark that the verification
of the $L_{2}-L_{2}$ bound on the operator $T$ is, in general, not trivial.
Intutively, $\theta (x, \lambda)$ constitute a generalized continuous
orhtogonal system, and hence the $L_{2}-L_{2}$ bound would be natural
if we had $\theta (x, \lambda),$ not $\theta (x, \lambda)^{2}$ in the definition
of the operator $T.$ In this respect, the problem of the perturbation of
the free Schr\"odinger operator has ``special status": indeed, we can take
in this case $\theta (x, \lambda)=\exp(i\sqrt{\lambda}x).$ The free exponent
$\exp (i\sqrt{\lambda}x)$ has a unique property that its square is also a
free exponent and hence the $L_{2}-L_{2}$ bound for the operator $T$ comes
 for free in this case.
But in general, the square in (3)
makes the problem more delicate.

  The main new technique we develop to prove the criteria inculdes a
  condition ensuring the a.e. convergence of
the linear integral operators of rather general form. The result we prove seems to be new, although its formulation is
quite natural.
It generalizes some well-known
results on the a.e. convergence of Fourier integrals.

Applying the new criteria, we show that every potential $V(x)$, verifying
$|V(x)| \leq C(1+x)^{-\frac{2}{3}-\epsilon},$ preserves the absolutely continuous
spectrum of the free Hamiltonian, relaxing the condition given in \cite{Kis}.
The new technique naturally has a wider range of applications than
perturbations of the free Schr\"odinger operator. We also consider the perturbations
of the periodic Schr\"odinger operators and find that the new criteria may
be   applied in this case.   

We also show that all results we prove have natural analogs for the whole-axis
problem.

The paper is organized as follows: in the first section we formulate the
main results and recall the general framework of the approach: the consequences of the Gilbert-Pearson subordinacy theory and the Harris-Lutz
mehtod. In the second section we develop the basic new technique: a theorem on the a.e. convergence
of integral operators. In the third section, we complete the proofs of the abstract
criteria for the preservation of the absolutely continuous spectrum
and discuss applications.

\begin{center}
 \section{ Main results and basic technique}
\end{center}

\rm \normalsize

Let us denote by $d\rho^{H_{V}}$ the spectral measure associated with
the Schr\"odinger operator $H_{V}$ in a usual way (see, e.g., \cite{Lev}). We denote the components
in a standart decomposition of this measure into absolutely continuous,
singular continuous and pure point parts by $d\rho_{\rm{ac}}^{H_{V}},$ $d\rho_{\rm{sc}}^{H_{V}}$ and $d\rho_{\rm{pp}}^{H_{V}}$ respectively. The singular part of the spectral measure, which is a sum of two latter components, we denote by $d\rho_{\rm{s}}^{H_{V}}.$ We use notation $m(E)$ and $\chi (E)$
for the Lebesgue measure and the characteristic function of the measurable set $E.$ One of our goals
in this paper is to prove the following theorem.

\noindent  \bf Theorem 1.1. \it Suppose that the potential $V(x)$
satisfies $|V(x)| \leq C(1+x)^{-\frac{2}{3}-\epsilon},$ where $\epsilon$ is
an arbitrary small positive number. Then the absolutely continuous part, $\rho_{\rm{ac}}^{H_{V}},$ of the spectral measure of the operator $H_{V},$ fills the
whole positive semi-axis \rm ( \it i.e. we have $\rho_{\rm{ac}}^{H_{V}}(T)>0$ for every
measurable $T \subset R^{+}$ with $m(T)>0$ \rm ).

\rm The second theorem  we prove here has to do with the slowly decreasing
perturbations of the Schr\"odinger operators with periodic potentials.
Suppose that $U(x)$ is a piecewise continuous periodic function on the real
axis and $H_{U}$ is defined, as before, on the half-axis. It is a well-known
fact that the spectrum of $H_{U}$ consists of the bands $\{ [a_{n}, b_{n}] \}^{\infty}_{n=0},$ where the spectrum is simple and purely absolutely continuous
and perhaps single eigenvalues in the gaps of the spectrum. This follows easily, for example, from the considerations in \cite{ReSi}. We have

\noindent \bf Theorem 1.2. \it Let $U(x)$ be piecewise continuous and periodic and let the set $S = \bigcup^{\infty}_{n=0} [a_{n}, b_{n}]$ be the absolutely continuous part of the spectrum of the operator $H_{U}.$ Then the absolutely continuous spectrum of $H_{U}$ is stable under all perturbations $V(x)$ satisfying
$|V(x)| \leq C(1+x)^{-\frac{2}{3}-\epsilon}$ (i.e. $\rho_{\rm{ac}}^{H_{U+V}}(T)>0$ for every $T \subset S$ with $m(T)>0$).  \\

\rm Of course, it may be more natural to consider an analog of Theorem 1.2 for the whole line problem. This poses no difficulty, as well as an extension of Theorem 1.1 to this case. Let us denote the Schr\"odinger operator given by the
differential expression $-\frac{d^{2}}{dx^{2}} + V(x)$ on the whole axis by $\tilde{H}_{V}.$ Then the following statements hold true: \\

\noindent \bf Theorem 1.3. \it Suppose that the potential $V(x)$ satisfies
$|V(x)| \leq C(1+x)^{-\frac{2}{3}-\epsilon}.$ Then the absolutely continuous
spectrum of the operator $\tilde{H}_{V}$ fills the whole positive semi-axis
with multiplicity two.  \\

\noindent \bf Theorem 1.4. \it Let $U(x)$ be piecewise continuous and periodic
and let $S=\bigcup^{\infty}_{n=0} [a_{n}, b_{n}]$ be the spectrum of the operator
$\tilde{H}_{U}.$ Then for every perturbation $V(x),$ satisfying $|V(x)| \leq
C(1+x)^{-\frac{2}{3}-\epsilon},$ the absolutely continuous spectrum of the operator $\tilde{H}_{U+V}$ fills $S$ with multiplicity two. \\

\rm These theorems will follow from the general criteria for the preservation of the absolutely continuous spectrum under slowly decreasing perturbations. We will give a precise formulation of these criteria in the end of this section when all relevant notation will be introduced in
the cousre of work and some of the technique will be developed which makes the
results more transparent.

The first important ingredient of the approach to the general slowly decaying
perturbations of the Schr\"odinger operators which we develop here is the relation
between the behavoir of the generalized eigenfunctions and the spectral
properties of the Schr\"odinger operators. Namely, we use the fact that for
a large class of potentials to show that a certain set $S$ belongs to the essential support of the absolutely continuous part of the spectrum of the 
operator $H_{W}$ it
suffices to show that for every energy from the set $S,$ all solutions of the
equation
\begin{equation}
 \left(-\frac{d^{2}}{dx^{2}} + W(x)\right)u = \lambda u
\end{equation}
 are bounded.
This result was first proven by Stolz \cite{Stolz} for the potentials $W(x)$ satisfying
\begin{equation}
W \in L^{1, \rm{loc}} \:\,\, \rm{and} \:\,\, \sup_{x} |\int_{|x-y| \leq 1} W_{-}(y)\,dy|< \infty,
\end{equation}
 where $W_{-}$ is the negative part of the potential $W.$ The proof of  Stolz relies on the Gilbert-Pearson subordinacy theory \cite{GiPe}, a rather recent development in the spectral theory of one-dimensional Schr\"odinger operators, which is remarkable in a way that it provides a direct and efficient relation between the properties of the solutions of (4) and the spectrum. The results
of \cite{GiPe} were elaborated and proofs simplified by Jitomirskaya and
Last in \cite{JiLa}. We will need the following \\

\noindent \bf Lemma 1.5. \it Suppose that the potential $W$ satisfies \rm (5)
\it and for every energy $\lambda$ from a certain set $S$ all solutions of the
 equation \rm (4) \it are bounded as $x \rightarrow \infty.$ Then the absolutely continuous spectrum of the operator $H_{W}$ fills the whole set $S,$ so that we have $\rho_{\rm{ac}}^{H_{W}}(T)>0$
for every measurable $T \subset S$ with $m(T)>0$. Moreover, no part of
the singular measure is supported on $S:$ $\rho_{\rm{s}}^{H_{V}}(S)=0.$ \\

\noindent \rm We remark that a simplest proof of this statement (with slightly more restrictive conditions on the potential) may be found in \cite{Sim3}. From the considerations in \cite{Sim3} also follows the analogous result for the whole-line problem: \\

\noindent \bf Lemma 1.6. \it Let the potential $W$ satisfy \rm (5). \it
Suppose that for every energy from the sets
$S_{\pm}$ all solutions of the generalized eigenfunction equation \rm (4)
\it are bounded as $x \rightarrow \pm \infty$ respectively. Then the absolutely
continuous spectrum of the operator $\tilde{H}_{W}$ fills the set $S_{+} \cup
S_{-}$
with the multiplicity at least one and singular component of the
spectral measure gives zero weight to the set $S_{+}\cup S_{-}:$ $\rho_{s}(S_{+}
\cup S_{-})=0.$
 Moreover, the absolutely continuous spectrum of the operator
$\tilde{H}_{W}$ fills the set $S_{+}\cap S_{-}$ with multiplicity two.  \\

\rm Now let us now formulate an additional assumption on the absolutely continuous
spectrum of the background operators $H_{U}$ and $\tilde{H}_{U}.$ \\

\noindent \bf Assumption. \it For the half-axis problem, we suppose that there exists
a measurable set $S$ of positive Lebesgue measure
such that for every energy $\lambda \in S$ all solutions of the equation \rm
(4) \it are bounded.

For the whole axis problem, we assume that there exist two sets $S_{\pm}$
such that
all solutions of the equation \rm (4) \it are bounded as $x \rightarrow \pm \infty$
for $\lambda \in S_{\pm}$ respectively, and $m(S_{+}\cup S_{-})>0.$ \\

\rm \noindent Notice that Lemmas 1.5 and 1.6 imply that the sets $S$ (respectively
$S_{\pm}$) belong to the
support of the absolutely continuous part of the spectral measure of the
operator $H_{U}$ (respectively $\tilde{H}_{U}).$ Moreover, the spectrum
on these sets is purely absolutely continuous, which means that the singular
parts of the corresponding spectral measures give no weight to these sets.
Here we will study the stability under slowly decaying perturbations
of the absolutely continuous spectrum of the operator $H_{U}$
(or $\tilde{H}_{U}$) supported on
$S$ ($S_{+}\cup S_{-}$) respectively. In other words, the methods we develop here are applicable only to the
type of the absolutely continuous spectrum which corresponds to the situation
when all solutions for the corresponding energies are bounded (at least on one of
the half-axes in the whole-line problem case). This condition
is not very restrictive. In most situations when one-dimensional
Schr\"odinger operators are known to have absolutely continuous spectrum,
it is exactly of the type described by the assumption. The question whether
there exist at all Schr\"odinger operators with the absolutely continuous spectrum of the different
type, i.e. such that for a set of the energies of positive measure from
the esential support of $\rho_{\rm{ac}}^{H_{U}}$ there exist unbounded solutions of
(4) was open for a long time. Although it  was  settled in a positive
way in \cite{Mol}, the corresponding examples are rather special
(for instance, the potential is unbounded both from above and from below).

\rm Now suppose that $H_{U}$ is a Schr\"odinger operator with 
part of the absolutely continuous spectrum supported on the set $S$ as in assumption.
Lemmas 1.5 and 1.6 reduce the problem of spectral analysis of the absolutely
continuous spectrum of $H_{U+V}$ to studing the asymptotics of the solutions
of the equation 
\begin{equation}
-\phi'' +(U+V)\phi =\lambda \phi.
\end{equation}
Namely, if we could show that all solutions of this equation are still bounded
for $\lambda \in S,$ we could apply again Lemma 1.5 to infer that the absolutely
continuous spectrum on $S$ is preserved. It is easy to show that for the short
range potentials this idea may be realized. However, as we mentioned in the introduction,
already for potentials decaying at a Coulomb rate imbedded singular spectrum
may appear, and it may be dense already for potentials decaying arbitrarily
slower than Coulomb. Hence the boundedness of all solutions may not hold already
in this case for rather rich set.
The main idea is now to show that for \it almost every \rm energy $\lambda$
from the set $S,$ all solutions of the equation (6) are bounded.

\rm  The second essential component of the approach is a certain asymptotic
integration method. We use it in a form proposed originally by Harris and
Lutz \cite{HaLu} for studing the asymptotics of the solutions of Schr\"odinger
equation with some particular oscillating potentials.
 We rewrite the
equation (6) as a system
\[ y' = \left( \begin{array}{cc} 0 & 1 \\ U+V-\lambda & 0 \end{array} \right)y. \]
We apply a variation of the parameters transformation
\[ y = \left( \begin{array}{cc} \theta (x, \lambda) & \overline{\theta} (x, \lambda)
\\ \theta '(x, \lambda) & \overline{\theta}' (x, \lambda) \end{array} \right) z, \]
to bring the system to the more symmetric form
\begin{equation} z' = \frac{1}{\Im (\theta \overline{\theta}')}\left( \begin{array}{cc} V(x)|\theta (x, \lambda)|^{2} & V(x) \overline{\theta} (x, \lambda)^{2} \\ -V(x) \theta (x, \lambda)^{2} & -V(x) |\theta (x, \lambda)|^{2} \end{array} \right)z. \end{equation}
Let us introduce a short-hand notation for the functions appearing in a latter
system, namely, let us write the system as
\begin{equation} z' = \left( \begin{array}{cc} D & L \\ \overline{L} & \overline{D} \end{array} \right)z,
\end{equation}
where $D= \frac{1}{\Im (\theta \overline{\theta'})}
V(x)|\theta(x, \lambda)|^{2}$ and
$L(x) = \frac{1}{\Im (\theta \overline{\theta})} V(x) \overline{\theta} (x, \lambda)^{2}.$
 The main approach to the study of the asymptotics of solutions for
systems similar to (7) is to attempt to find some transformation which will
reduce the off-diagonal terms so that they will become absolutely integrable
and then try to apply Levinson's theorem \cite{Lev} on the $L^{1}$-perturbations
of the systems of linear differential equations. Indeed, if we could assume that
the off-diagonal terms of the system (7) are absolutely integrable, then in many cases the main term of the asymptotics of the solutions of the system (7) would, by
Levinson's theorem (see, e.g. \cite{Lev}), coincide with the
solution of the system (7) with only diagonal part present. These solutions
are all bounded and going back to the original Schr\"odinger equation we also find the asymptotics of its solutions and see that they are also bounded.

It was discovered by Harris and Lutz \cite{HaLu} that when $W(x)$ is
a conditionally integrable function, the following simple transformation
of the system (8) works in some cases. We let
\begin{equation}
z(x)= (1-|q|^{2})^{-\frac{1}{2}}(I + {\cal Q})\omega(x)
\end{equation}
where $I$ is an identity matrix, while ${\cal Q}$ is given by
\[ {\cal Q}(x)= \left( \begin{array}{cc} 0 & q(x) \\ \overline{q}(x) & 0
\end{array} \right) \]
where $q(x)= -\int^{\infty}_{x} L(y)\,dy.$
In this case $q(x) \stackrel{x \rightarrow \infty}{\longrightarrow} 0$,
so that for large enough $x$ the transformation (9) is non-singular
and preserves the asymptotics of the solutions. For the
new variable $\omega(x)$ we obtain a system
\begin{equation}
\omega ' = \left( \left( \begin{array}{cc} D & 0 \\ 0 & \overline{D} \end{array}
 \right) + (1-|q|^{2})^{-1}) \left( \begin{array}{cc} i \Im (q \overline{L})+
2|q|^{2}D & 2\overline{qD} - \overline{q}^{2}L \\ 2qD-q^{2}\overline{L}
 & -i \Im (q \overline{L}) +  2|q|^{2} \overline{D} \end{array} \right) \right) \omega.
\end{equation}
Because of the way the function $q$ was defined, we see that the terms in the
second summand, which we would like to treat as a perturbation, decay faster
than the elements of the original system. In particular, for the system (10)
that we consider, for the energies for which we manage to define the
function $q(x, \lambda)$ as above, every entry of the second matrix 
is a product of the function
\[  V(x)
\int\limits^{\infty}_{x} V(t) \theta (x, \lambda)^{2}\,dt =
V(x) q(x, \lambda) \]
and some bounded function of $x.$

It was shown in \cite{Kis} that in the case of the background potential
$U$ equal to zero if we take $\theta (x, \lambda)=\exp (i\sqrt{\lambda})x$ 
and the potential $V(x)$ satisfies
$|V(x)| \leq C(1+x)^{-\frac{3}{4}-\epsilon},$ then for a.e.~$\lambda \in (0, \infty)$
we have $q(x, \lambda)V(x) \in L^{1}.$ Hence, for these $\lambda$
we can apply the Levinson's theorem  (or, in our particular case,
more straightforward integral equation technique as in
\cite{ReSi3}), to find the asymptotics of the solutions and see that all
solutions are bounded for a.e. positive $\lambda.$ 
The problem of proving that $V(x)q(x, \lambda) \in L^{1}$ reduces to 
studing the function $q(x, \lambda),$ in particular the a.e.~$\lambda$ 
existence
and rate of the convergence of the integral defining $q(x, \lambda).$
We remark that when $U=0,$ the function $q(x, \lambda)$ is just the ``tail" of the 
Fourier transform of the potential $V(x).$ Parseval formula may be used
in this case not only to show that $V(x)q(x, \lambda) \in L^{1}$
for a.e.~$\lambda,$ but also to describe rather explicitely in terms
of the Fourier transform of $x^{\frac{1}{4}}V(x)$ the exceptional 
divergence set where the singular part of the spectral measure might 
be supported.

If we study the structure of the perturbation in the system (10), we
see that the slowest decaying off-diagonal terms are of type $2qD$ or
$2\overline{q}\overline{D}.$ If we had $D=0,$ the slowest decaying terms
would be $q^{2}\overline{L}$ and $\overline{q}^{2}L,$ which contain the
decaying function $q$ in the power two. Now let us perform one more
transformation with the system (8):
\[ z= \left( \begin{array}{cc} \exp \left( -\frac{1}{\Im (\theta \overline{\theta}')}
\int^{x}_{0}V(t)|\theta(t, \lambda)|^{2}\,dt \right)  & 0 \\ 0 & \exp \left( \frac{1}{\Im (\theta \overline{\theta}')} \int^{x}_{0} V(t)|\theta(t, \lambda)|^{2}\,dt \right) \end{array} \right)z_{1} . \]
Let us denote by $L_{1}(x, \lambda)$ the function
\begin{equation} 
\frac{1}{\Im (\theta \overline {\theta}')} V(x) \overline{\theta} (x, \lambda)^{2} \exp \left( - \frac{2}{\Im (\theta \overline{\theta}')} \int^{x}_{0}
V(t)|\theta(t, \lambda)|^{2}\,dt \right).
\end{equation}
For the variable $z_{1}$ we have
\[ z_{1}' =  \left( \begin{array}{cc} 0 & L_{1}(x, \lambda) \\ 
\overline{L}_{1}(x, \lambda) & 0 \end{array} \right)z_{1}. \]
Define
\begin{equation}
 q_{1}(x, \lambda) = - \int\limits_{x}^{\infty}L_{1}(s, \lambda)
\,ds.
\end{equation}
(for those $\lambda$ for which such integral will exist) and perform
an $I+Q_{1}$ transformation
\[ z_{1}  = (1-|q_{1}|^{2})^{-\frac{1}{2}} \left( \begin{array}{cc}
1 & q_{1} \\ \overline{q}_{1} & 1 \end{array} \right) \omega_{1}. \]
We get for the new variable $\omega_{1}$
\begin{equation}
\omega_{1}' = (1-|q_{1}|^{2})^{-1} \left( \left( \begin{array}{cc}
\frac{1}{2} (q_{1} \overline{L}_{1} - \overline{q}_{1}L_{1}) & 0 \\
0 & -\frac{1}{2}(q_{1}\overline{L}_{1} - \overline{q}_{1}L_{1})
\end{array} \right) + \left( \begin{array}{cc} 0 & \overline{q}_{1}^{2}L_{1} \\
q_{1}^{2}\overline{L}_{1} & 0 \end{array} \right) \right) \omega_{1}.
\end{equation}
Here the first diagonal term has purely imaginary entries and hence leads
to bounded solutions, while the second term decays as $|q_{1}(x, \lambda)|^{2}V(x).$
This computation suggests that we should study the question of existence
and convergence of the integrals like
\[ q(x, \lambda) = \int\limits_{x}^{\infty}V(t) \overline{\theta} (t, \lambda)^{2}\,dt \]
and 
\[ q_{1}(x, \lambda) = \int\limits_{x}^{\infty}V(t)\overline{\theta}(t, \lambda)^{2} \exp \left( - \frac{2}{\Im (\theta \overline{\theta}')} \int\limits^{x}_{0}
V(t)|\theta(t, \lambda)|^{2}\,dt \right). \]
In particular if we could get the a.e. convergence estimates for latter integrals
similar to the case of Fourier integral, then we could 
study the asymptotics of solutions of Schr\"odinger operators with
potentials decaying slower than at $x^{-\frac{3}{4}-\epsilon}$ rate.

The next section is devoted to the handling of the a.e. convergence questions
and reduction to a simpler issue of norm estimates.

We conclude this section with the formulation of the two general criteria 
for the preservation of the absolutely continuous spectrum, both for half-line and full line problems.

 Let $H_{U},$ $\theta (x, \lambda)$ and $S$ be as above. 
It will be convenient for us to choose a certain basis of solutions for every
$\lambda \in S.$ Namely, for any $\lambda \in S$ we choose two linearly
independent, measurable in $(x, \lambda)$ solutions $\theta (x, \lambda)$ and its complex conjugate
$\overline{\theta}(x, \lambda)$ such that:
\begin{equation}
|\theta(x, \lambda)| \leq C \,\, \rm{for} \,\, \rm{every} \,\, \lambda \,\,  \rm{uniformly}
\,\,in \,\, \lambda \in S.
\end{equation}
 Of course we can always choose such basis $\theta(x, \lambda),$
$\overline{\theta}(x, \lambda)$ if all solutions of (4) are bounded for every
$\lambda \in S.$

For the whole-axis problem, for every $\lambda \in S_{\pm}$ we choose a
basis consisting of solutions $\theta_{\pm}(x,\lambda),$ $\overline{\theta_{\pm}}
(x, \lambda)$ satisfying (14) when $x \in (0, \pm \infty)$ respectively.

For example, treating the free case we will have $S=(0,\infty)$ and 
$\theta (x, \lambda) = \exp (i \sqrt{\lambda}x)$ and considering the periodic
case we will have $S=\cup_{n=1}^{\infty} (a_{n}, b_{n})$ and $\theta (x, \lambda)$ will be chosen to be Bloch functions.

We now perturb the
operator $H_{U}$ by a decaying potential $V(x).$  Define the linear
operators $T_{1}$
and $T_{2}$ acting
on the bounded functions of compact support by 
\begin{equation} (T_{1}f)(\lambda) = \frac{\chi (S)}{\Im (\theta \overline{\theta}')} \int\limits_{0}^{\infty}\overline{\theta}(x, \lambda)^{2}f(x)\,dx
\end{equation}
and 
\begin{equation}
(T_{2}f)(\lambda)=  \frac{\chi (S)}{\Im (\theta \overline{\theta}')} \int\limits_{0}^{\infty}\overline{\theta}(x, \lambda)^{2}
\exp  \left( - \frac{2}{\Im (\theta \overline{\theta}')} \int\limits^{x}_{0}
V(t)|\theta(t, \lambda)|^{2}\,dt \right) f(x)\,dx.
\end{equation}
We have \\

\noindent \bf Theorem 1.7. \it Suppose that there exists the partition of the set $S,$ $S=\cup_{i=1}^{\infty} S_{i},$ such that for every $i$ we have the bounds
\[ \|\chi(S_{i})T_{1}f\|_{2} \leq C_{1i}\|f\|_{2}.\]
Then the absolutely continuous spectrum of the operator $H_{U},$ supported
on $S,$ is stable under perturbations by all potentials $V(x)$ satisfying
$|V(x)| \leq C(1+x)^{-\frac{3}{4}-\epsilon}$ \rm ( \it i.e. $\rho_{\rm{ac}}^{H_{U+V}}
(S_{1})>0$ for every $S_{1} \subset S$ with $m(S_{1})>0$ \rm ). \it Moreover, for a.e.~$\lambda \in S$ we have two linearly independent solutions
$\phi (x, \lambda),$ $\overline{\phi} (x, \lambda)$ of the equation 
$(H_{U+V} -\lambda) \phi = 0$ with the asymptotics
\begin{equation}
 \phi (x, \lambda) = \theta (x, \lambda) \exp  \left( - \frac{2}{\Im (\theta \overline{\theta}')} \int\limits^{x}_{0}
V(t)|\theta(t, \lambda)|^{2}\,dt \right)\left(1+ O(x^{-\epsilon})\right). 
\end{equation}

\noindent \bf Theorem 1.8. \it Assume that a potential $V(x),$
verifies $|V(x)| \leq C(1+x)^{-\frac{2}{3}-\epsilon}.$ Suppose that there exists
a partition of the set $S$ into the sets $\{ S_{i} \}_{i=1}^{\infty}$ such that for every $i$ we have the bounds
\[ \|\chi(S_{i})T_{2}f\|_{2} \leq C_{2i}\|f\|_{2}. \]
Then the set $S$ remains in the support of the absolutely continuous 
spectrum of the operator $H_{U+V}$ \rm ( \it i.e. $\rho_{\rm{ac}}^{H_{U+V}}(S_{1})>0$
for every measurable $S_{1} \subset S$ with $m(S_{1})>0$ \rm ). \it Moreover, 
for a.e.~$\lambda \in S$ we have two linearly independent solutions 
$\phi (x, \lambda),$ $\overline{\phi} (x, \lambda)$ of the equation 
$(H_{U+V} -\lambda) \phi = 0$ with the asymptotics
\begin{eqnarray}
 \phi (x, \lambda) & = & \theta (x, \lambda) \exp  \left( - \frac{2}{\Im (\theta \overline{\theta}')} \int^{x}_{0}
V(t)|\theta(t, \lambda)|^{2}\,dt + i\int\limits_{0}^{x}\Im(q_{1}(t, \lambda)\overline{L}_{1}(t, \lambda))\,dt  \right) \times \nonumber \\ 
& & \times \left(1+ O(x^{-\epsilon})\right),
\end{eqnarray}
where $L_{1}$ and $q_{1}$ are as in \rm (11) \it and \rm (12). \\

The whole-axis analogs of Theorems 1.7 and 1.8 are formulated 
as follows. Let, as before, $\tilde{H}_{U}$ be a Schr\"odinger operator
with potential $U$ defined on the whole axis. Let $S_{-},$ $S_{+},$
$\theta_{-}(x, \lambda),$ and $\theta_{+}(x, \lambda)$ be as above. 
 Let the operators $T^{\pm}_{1},$ $T^{\pm}_{2}$ be defined on bounded functions of compact
support by
\[ (T^{-}_{1}f)(\lambda) = \frac{\chi (S_{-})}{W[\theta_{-}, \overline{\theta}_{-}]}
\int\limits_{-\infty}^{0}\overline{\theta}_{-}(x, \lambda)^{2}f(x)\,dx \]
and
\[ (T_{2}^{-}f)(\lambda)=  \frac{\chi (S_{-})}{\Im (\theta_{-} \overline{\theta}_{-}')} \int\limits_{-\infty}^{0}\overline{\theta}_{-}(x, \lambda)^{2}
\exp  \left( - \frac{2}{\Im (\theta_{-} \overline{\theta}_{-}')} \int\limits^{x}_{0}
V(t)|\theta_{-}(t, \lambda)|^{2}\,dt \right) f(x)\,dx . \]
In the definitions of $T_{1,2}^{+}$ we just replace all signs ``$-$"
by the ``$+$" signs in the right-hand side.

Then we have the following criteria: \\

\noindent \bf Theorem 1.9. \it Suppose that there exist partitions of the sets $S_{\pm}$ into countable unioins of sets $\{ S_{\pm i} \}_{i=1}^{\infty}$ respectively such that the following bounds hold for every $i:$
$\| \chi (S_{+ i})T^{+}_{1}f\|_{2} \leq C_{+i}\|f\|_{2}$ and $\|\chi(S_{- i}) T^{-}_{1}f\| \leq C_{-i}\|f\|_{2}$ for every $i.$ Then the absolutely continuous
spectrum of the operator $H_{U},$ supported on $S_{-}$ and $S_{+},$ is stable under perturbations
by all potentials $V(x)$ verifiyng $|V(x)| \leq C(1+|x|)^{-\frac{3}{4} -\epsilon}.$ 
Namely, for every such $V$ the operator $\tilde{H}_{U+V}$ has absolutely
continuous spectrum of multiplicity at least one on  
$S_{-} \cup S_{+}$ and of multiplicity two on $S_{+} \cap S_{-}.$ \\

\noindent \bf Theorem 1.10. \it Let potential $V(x)$ satisfy $|V(x)| \leq C(1+|x|)^{-\frac{2}{3}-\epsilon}.$ Suppose that there exist the partitions
of the sets $S_{\pm}$ into the countable unions of the sets $\{ S_{\pm i}\}_{i=1}^{\infty}$ respectively so that the following bounds hold for every
$i:$
$\| \chi(S_{+ i})T^{+}_{2}f\|_{2} \leq C_{+i}\|f\|_{2}$ and $\|\chi (S_{- i})T^{-i}_{2}f\| \leq C_{-i}\|f\|_{2}$ for every $i.$ Then the absolutely continuous
spectrum of the operator $\tilde{H}_{U},$ supported on $S_{-}$ and $S_{+},$ is preserved under pertubation by the potential $V(x).$
Namely, the operator $\tilde{H}_{U+V}$ has absolutely
continuous spectrum of multiplicity at least one on $S_{+} \cup S_{-}$  and of multiplicity two on $S_{+} \cap S_{-}.$ \\

\rm All these criteria show that instead of studying the rate of a.e. convergence
of certian integral operators all we have to do is to check a certain 
$L_{2}-L_{2}$ estimate, which is in many situations much simpler. The passage from a.e. convergence questions to norm estimates
is the third ingredient of our approach. This is a crucial element and the main
new idea in the context of the spectral study of Schr\"odinger operators.
We treat this subject in the next section.
 
\begin{center}
 \section{ A.e. convergence for integral operators}
\end{center}

\normalsize

Let the operator $T$ be defined on the measurable bounded functions $f$
of compact support by
\begin{equation}
(Tf)(k) = \int\limits_{0}^{\infty} A(k,x) f(x)\, dx,
\end{equation}
where $A(k,x)$ is a measurable and bounded function on $R^{2}.$
Let us denote by $A$ the upper bound on the kernel $A(k,x).$
Denote by $Mf(k)$ the corresponding maximal function
\begin{equation}
Mf(k) = \sup_{N} \left| \int\limits_{0}^{N}A(k,x)f(x)\,dx \right|.
\end{equation}
We are interested in studying the a.e. convergence questions for the integral
operators of type (19),
i.e. in finding more or less simple conditions implying that
the integral (19) converges for a.e. $k$ if, say, $f$ belongs to some $L^{p}$ space. The natural class of conditions to look at are norm estimates on the
operator $T.$ When we deal with Fourier integrals, in many instances a Parseval
equality with explicit kernel helps us to get some information about a.e.
convergence (see, e.g., \cite{Tit}).
Of course, there is rarely such thing as the Parseval equality for the integrals
of type (19). However, a weaker tool - norm estimates - turns out to be 
sufficient for our needs. We use a traditional method of maximal function
estimates to infer the a.e. convergence. The principal result of this
section is the derivation of the estimates on the maximal function (20) given certian norm estimates on the operator
itself. This reduces the proof of the a.e. convergence of the operator $T$
given by (19) for the functions from the certain class to establishing the appropriate
norm estimates for $T.$

To formulate the main result of this section in the natural form, it is
useful to introduce the scale of Lorentz spaces $L_{pq}.$ We remind here the basic
definitions and properties of these spaces. For more details and proofs
we refer to \cite{Stein}. The function $f$ belongs to $L_{pq}$ iff
\[ \|f\|^{*}_{pq} = \left( \frac{q}{p}\int\limits^{\infty}_{0} [t^{\frac{1}{p}}
f^{*}(t)]^{q}\frac{dt}{t} \right)^{\frac{1}{q}} < \infty , \]
where $f^{*}(t)$ is a non-decreasing rearrangement of the function $f,$
defined by
\[ f^{*}(t) = \inf \{ s: m\{x| \, |f(x)|>s \} \leq t \}. \]
The expression $\|f\|^{*}_{pq}$ does not generally define a norm since
it does not generally satisfy the triangle inequality. However, if $1 \leq p
\leq \infty$ and $1<q\leq \infty,$ there exists
a norm $\|\,\,\|_{pq}$ on $L_{pq}$ which is equivalent to $\|\,\,\|^{*}_{pq}.$
In particular, $\|f\|^{*}_{pq} \leq \|f\|_{pq} \leq \frac{p}{p-1} \|f\|^{*}_{pq},$
see \cite{Stein}.
The $L_{pq}$ spaces are of ``$L_{p}$- type", i.e. for every characteristic
function of the measurable set $E$, $\chi (E),$ we have $\| \chi (E) \|_{pq}=
(m(E))^{\frac{1}{p}}.$ Also $\| f \|^{*}_{pq_{1}} \leq \|f\|^{*}_{pq_{2}}$ if $q_{1}
>q_{2},$ so that for fixed $p$, the Lorentz space extends as $q$ grows. Finally, we note that the $L_{pp}$-space coincides with the
usual $L_{p}$ since $\|f\|_{pp} = \|f\|_{p}.$

  We will also need the following well-known interpolation theorem, so-called
  generalized Marcinkiewicz theorem. We refer
to \cite{Stein} or \cite{BeL} for a proof. The operator $B$ is called subadditive if it
satisfies $|B(f_{1}+f_{2})(k)| \leq |Bf_{1}(k)| +|Bf_{2}(k)|.$
We say that $B$ is of restricted weak type $(r,p)$ if its domain $D(B)$ contains
all finite linear combinations of characteristic functions of sets of finite
measure and all truncations of its members and satisfies $\|Tf\|_{p\infty}
\leq C \|f\|_{r1}$ for all $f \in D \cap L_{r1}.$ \\

\noindent \bf Theorem. \it Suppose that $T$ is a subadditive operator
of restricted weak types $(r_{j}, p_{j}),$ $j=1,2$ with $r_{0}<r_{1}$ and
$p_{0} \neq p_{1},$ then there exists constant $B_{\theta}$ such that
\[ \| Tf \|_{pq} \leq B_{\theta} \|f \|_{rq} \]
for all $f$ from the domain of $T$ and $L_{rq},$ where $1 \leq q \leq \infty,$
$\frac{1}{p} = \frac{1-\theta}{p_{0}}+ \frac{\theta}{p_{1}},$ $\frac{1}{r} =
\frac{1-\theta}{r_{0}}+\frac{\theta}{r_{1}},$ $0<\theta<1.$ \\

\rm Now we are in a position to formulate the main result that we need. The theorem
below generalizes the classical theorem of Zygmund \cite{Zyg} for the case of Fourier
transform. 

\noindent \bf Theorem 2.1. \it Suppose that an operator T, defined by \rm (19) \it with bounded kernel $A(k,x),$ satisfies the norm estimate
$\| Tf\|_{2} \leq C_{q} \|f\|_{2}$ for all bounded functions of compact
support. Then for every $q>2$ and $p$ such that $q^{-1}+p^{-1}=1,$
we have the following estimate for the maximal function:
\begin{equation}
 \|Mf\|_{q} \leq C_{q} \|f\|_{pq} \;\: for \;\: every \;\: f \in L_{pq}
\end{equation}
{\rm (} hence, in particular, $\|Mf\|_{q} \leq C_{q}\|f\|_{p}${\rm )}.
As a consequence, the integral 
\[ \int^{N}_{0} A(k,x)f(x)\,dx \]
 converges
as $N \rightarrow \infty$ for almost every value of $k$ if $f \in L_{pq}$
{\rm (} and in particular if $f \in L_{p},$ $1 \leq p < 2${\rm )}. \\

\noindent \it Remark. \rm The result we prove here is suited for
the applications we make in this paper. In fact, similar results hold
in greater generality for more general integral operators with natural defintions of the
maximal function. The proof is  more involved and we plan to
devote a separate publication \cite{Kis1} to this problem.

\noindent \bf Proof. \rm First we remark that from the estimate on the
maximal function the a.e. convergence for operator $T$ in classes $L_{pq},$ $p<2,$ follows in a standart way (see,
e.g., \cite{Gar}). We sketch here this simple argument for the sake of completness. Indeed,
suppose that there exists a function $f \in L_{pq}$ with $p$ and $q$ as
in the theorem and a set $S$ of positive measure where the integral defining
$Tf(k)$ diverges. Then we can find $\epsilon >0$ and a set $S_{\epsilon}
\subset S,$ such that $m(S_{\epsilon})>0,$ and for every $k \in
S_{\epsilon}$ and every positive number $N_{1}$ there exists a
larger number $N_{2}$ such that $|\int_{N_{1}}^{N_{2}}A(k,x)f(x)\,dx|>
\epsilon.$ Hence, for every $N_{1},$ we have $\|Mf \chi (N_{1}, \infty) \|
_{q} \geq \epsilon m(S_{\epsilon})^{\frac{1}{q}}.$ On the other hand,
clearly $\|f\chi (N_{1}, \infty) \|_{pq} \stackrel{N_{1}
\rightarrow \infty}{\longrightarrow} 0.$ This contradicts the bound (21).

We now come to a proof of the first assertion. The first step is the  decompostion of the support of the function $f$ into dyadic pieces and estimates on the certain auxilliary  maximal functions. A similar idea was  used already by Paley \cite{Pal} in his proof of a.e. convergence for the series of orhtogonal functions. Let $f$ be a measurable bounded function of compact
support and choose $n$ so that $2^{n-1} \leq m(\rm{supp}(f)) \leq 2^{n}.$ 
 Let the measurable set $E$ be the
support of the function $f:$ $E= \{ x| \;\, |f(x)|>0 \}.$ For every
integer $m<n,$ we consider a partition of the set $E$ into the following
sets $E_{m,l}:$
\[ E_{m,l} = (a_{m,l}, a_{m,l+1}) \cap E, \]
where $a_{m,l}$ is defined by a condition
\[ a_{m,l} = \inf \{ a | m((0,a) \cap E) \} = 2^{m}l. \]
The number of the sets $E_{m,l}$ is between $2^{n-m}$ and $2^{n-m-1}.$
For notational convenience, we will assume that this number is always
$2^{n-m}$ and will define the missing $E_{m,l}$ as empty sets.

Let us define functions $M_{m,l}f$ and $M_{m}f$ by
\[ M_{m}f(k) = \sup_{l} |M_{m,l}f(k) | \]
and
\[ M_{m,l}f(k) = \left| \,\int\limits_{E_{m,l}} A(k,x)f(x)\,dx \right|. \]
Considering a dyadic development of every real number $N,$ it is easy to see
that
\begin{equation}
Mf(k) \leq \sum\limits^{n}_{m=-\infty} M_{m}f(k).
\end{equation}
Indeed, suppose that for a given value of $k,$ the supremum in (20) is reached
when the upper limit is $N$ (clearly if $f$ has compact support, the
supremum is reached for some value of $N$). Define the real number $s$ by
$s= m(E \cap (0,N))$ and consider the dyadic development
$s= \sum_{m=-\infty}^{n} s_{m}2^{m},$ where $s_{m}$ is equal
to $0$ or $1$ for every $m.$ Then by construction, we can find disjoint
sets $E_{m,l},$ at most one for each value of $m,$ so that $m((E\cap(0,N))/
(\cup_{m} E_{m,l}))=0.$ In fact, for each $m$ the set $E_{m,l}$ belongs to
the union iff $s_{m}=1.$ The corresponding value of $l$ then may be found
by the formula $l= \sum^{n}_{j=m+1} s_{j} 2^{j-m}.$

Obviously, we also have
\begin{equation}
 M_{m}^{q}f(k) \leq \sum\limits_{l=1}^{2^{n-m}} M_{m,l}^{q}f(k),
\end{equation}
for every $q>0.$ Fix now any $q>2$ and let $q'$ satisfy $q>q'>2.$
Under the conditions of the theorem, we have that $\|Tf\|_{\infty}\leq
C_{1}\|f\|_{1}$ for all $f \in L_{1}$ and $\|Tf\|_{2} \leq C_{2}\|f\|$ is
satisfied for all measurable bounded functions of compact support. By 
interpolation, we have
\[ \|Tf\|_{q'} \leq C_{q'} \|f\|_{p'q'}, \,\,\, (p')^{-1}+(q')^{-1}=1. \]
Noting that
\[ \|f\|^{*}_{p'q'} = \left( \frac{q'}{p'} \int\limits^{\infty}_{0} |f^{*}(t)|^{q'}
t^{\frac{q'}{p'}-1}\,dt\right) ^{\frac{1}{q'}}, \]
and using the equivalence of $\|\cdot \|^{*}$ and $\|\cdot\|,$ we see that in particular
\[ \|M_{m,l}f(k)\|^{q'}_{q'}= \|T(f \chi (E_{m,l}))(k) \|_{q'}^{q'} \leq
C_{q'}^{q'}\left( \frac{q'}{p'} \int\limits_{0}^{\infty} | (f\chi (E_{m,l}))^{*}
|^{q'}t^{\frac{q'}{p'}-1} \,dt \right)\leq \]
\[ \leq C_{q'}^{q'}2^{m(\frac{q'}{p'}-1)}\frac{q'}{p'} \| f\chi (E_{m,l})(t)\|^{q'}_{q'}.\]
Summing over $l$ and using (23), we obtain
\[ \|M_{m}f\|_{q'}^{q'} \leq C_{q'}^{q'}\frac{q'}{p'} 2^{m(\frac{q'}{p'}-1)}\|f\|_{q'}^{q'}. \]
By (22), we have that
\begin{equation}
\| Mf \|_{q'} \leq C_{q'} \left( \frac{q'}{p'} \right)^{\frac{1}{q'}} \|
f \|_{q'} \sum\limits_{m=-\infty}^{n} 2^{m (\frac{1}{p'}-\frac{1}{q'})}=
B_{q'} 2^{n(\frac{1}{p'}-\frac{1}{q'})}\|f\|_{q'}.
\end{equation}
Now we note that in a particular case when $f$ is a characteristic function
of a set, $f=\chi (E),$ (24) means
\begin{equation}
 \| M\chi (E) \|_{q'} \leq B_{q'} 2^{n(\frac{1}{p'}-\frac{1}{q'})}2^{\frac{n}{q'}}=
B_{q'} 2^{\frac{n}{p'}} \leq 2^{\frac{1}{p'}}B_{q'} \| \chi (E) \|_{p'1}. 
\end{equation}
It is easy to check that the operator $Mf,$ defined originally on the measurable
bounded functions of compact support, is a sublinear operator. 
It is  well-known fact that from the inequality (25) for sublinear operator it follows that
$\|Mf\|_{q'} \leq C\|\chi (E) \|^{*}_{p',1}$ holds for all finite combinations
of simple functions (see \cite{Stein}) and hence by simple limiting
argument for all measurable bounded functions of compact support. Interpolating with an
obvious relation $\|Mf\|_{\infty} \leq A \|f\|_{1},$ we obtain that $\|Mf\|_{\tilde{q}} \leq C_{\tilde{q}}\|f\|_{\tilde{p}
\tilde{q}}$ for every $\tilde{p},$ $\tilde{q}$ such that $q'>\tilde{q}>2$
and $\tilde{p}^{-1}+\tilde{q}^{-1}=1$ and for every function $f$ bounded
and of compact support. In particular, this relation holds for the value of $q$
we fixed in the beginning of the proof (and hence for every $q>2$).
 It is straightforward to see
that this inequality is then extended to all functions $f \in L_{pq}.$ $\Box$

\begin{center}
 \section{ Proofs of criteria and applications}
\end{center}

\normalsize

Now we prove all theorems formulated in the first section. First we give
proofs of the general criteria. \\

\noindent \bf Proof of Theorem 1.7. \rm By assumption, the operator $T_{1i}$
defined by 
\[ (T_{1i}f)(\lambda) = \frac{\chi (S_{i})}{\Im (\theta \overline{\theta}')}
\int\limits_{0}^{\infty} \overline{\theta}(x, \lambda)^{2} f(x)\,dx, \]
satisfies an $L_{2}-L_{2}$ bound on the bounded functions of compact support.
We may also assume that on every set $S_{i}$ we have the Wronskian $W[\theta, \overline{\theta}]$ bounded away from zero by some constant $c_{i},$ or else subdivide the $S_{i}$ so that it holds true.  Since the functions $\theta (x, \lambda)$ are uniformly bounded when $\lambda
\in S,$ we also have an obvious $L_{1}-L_{\infty}$ bound for $T_{1i}.$ 
Hence, by Theorem 2.1, the integral
 \begin{equation}
 \int\limits_{0}^{x} \overline{\theta}(y, \lambda)^{2} f(y)\,dy 
\end{equation}
converges as $x \rightarrow \infty$ for a.e.~$\lambda \in S_{i}$ for every $f \in L_{p},$ $1\leq p <2.$ By the 
assumption on the decay of $V(x),$ the function $x^{\frac{1}{4}}V(x)$ belongs
to $L_{2(1-\epsilon)}.$ Therefore, we find that 
\[ q(x, \lambda) = \int\limits_{x}^{\infty} \overline{\theta}(y, \lambda)^{2}
V(y)\,dy = \int\limits_{x}^{\infty} \overline{\theta}(y, \lambda)^{2} (V(y)
y^{\frac{1}{4}})y^{-\frac{1}{4}}\,dy= \]
\[ = x^{-\frac{1}{4}}\int\limits_{0}^{x} \overline{\theta}(y, \lambda)^{2}
(y^{\frac{1}{4}}V(y))\,dy + \frac{1}{4} \int\limits_{x}^{\infty} y^{-\frac{5}{4}}\int\limits_{0}^{y} \overline{\theta}(t, \lambda)^{2} (t^{\frac{1}{4}}V(t))\,dt. \]
We conclude that for all $\lambda \in S_{i}$ such that the integral (26) converges, and hence for a.e.~$\lambda \in S_{i},$ the function $q$ satisfies
$q(x, \lambda)=O(x^{-\frac{1}{4}})$ as $x \rightarrow \infty.$ Since this holds for any $i,$ we have that this estimate is also true for a.e.~$\lambda \in S.$ This implies
that $q(x, \lambda)V(x) \in L_{1}$ and allows for a.e.~$\lambda \in S$ to 
find the asymptotics of solutions of the perturbed Schr\"odinger equation.
Transforming back via (7) and (9) we find that for a.e.~$\lambda \in S,$
there exist two solutions $\phi(x, \lambda),$ $\overline{\phi}(x, \lambda)$
of the generalized eigenfunction equation (4) with the following asymptotics:
\[  \phi (x, \lambda) = \theta (x, \lambda) \exp  \left( - \frac{2}{\Im (\theta \overline{\theta}')} \int\limits^{x}_{0}
V(t)|\theta(t, \lambda)|^{2}\,dt \right)\left(1+ O(x^{-\epsilon})\right), \]
\[ \phi' (x, \lambda) = \theta' (x, \lambda) \exp  \left( - \frac{2}{\Im (\theta \overline{\theta}')} \int\limits^{x}_{0}
V(t)|\theta(t, \lambda)|^{2}\,dt \right)\left(1+ O(x^{-\epsilon})\right). \]
Clearly the solutions $\phi,$ $\overline{\phi}$ are linearly independent,
since the Wronskian $W[\phi, \overline{\phi}] = W[\theta, \overline{\theta}]
\ne 0.$ This concludes the proof, given Lemma 1.5. $\Box$ \\

\noindent \bf Proof of Theorem 1.8. \rm  Similarly to the previous proof,
we infer that under the assumption of the theorem, for every $f \in L_{p},$ $1 \leq p <2,$ the integral
\[ \int\limits_{0}^{x} \overline{\theta}(y, \lambda)^{2} \exp \left(-\frac{1}{\Im (\theta \overline{\theta}')} \int\limits_{0}^{y} V(t)
 |\theta(t, \lambda)|^{2}\,dt \right) f(y)\,dy \]
converges for a.e.~$\lambda \in S.$ As before, integrating by parts, we
find that if $|V(x)| \leq C(1+x)^{-\frac{2}{3}-\epsilon},$ then the
function $q_{1}(x, \lambda)$ given by (12) satisfies
\[ q_{1}(x, \lambda) = O(x^{-\frac{1}{6}}) \]
for a.e.~$\lambda \in S.$ Therefore, for a.e.~$\lambda \in S$ we also have
$|q_{1}(x, \lambda)^{2}V(x)| \leq C(1+x)^{-1-\epsilon}.$ This allows us to find
the asymptotics of the solutions of the system (13) and then of the original
Schr\"odinger equation. As in the previous proof, the Wronskian argument 
shows linear independence of the solution with asymptotics (18) and its
complex conjugate.  $\Box$ \\

The proofs of the whole line analogs of the criteria, Theorems 1.9 and 1.10, follow in the same way given Lemma 1.7. \\

Now we come to a final goal of this paper - concrete applications to the 
preservation of the absolutely continuous spectrum of the Schr\"odinger
operators. We first discuss the free case: $U(x)=0.$ We remark that the criterion
given by Theorem 1.7 applies trivially since the operator $T_{1}$ in question
is just a rescaled Fourier transform and hence satisfies the $L_{2}-L_{2}$ estimate.  This gives the stability of the absolutely continuous spectrum
of the free Schr\"odinger operators under perturbations by all potentials $V$
satisfying $V(x) \leq C(1+x)^{-\frac{3}{4}-\epsilon}.$ This has been proven 
in \cite{Kis} using a more direct method rather than Theorem 1.7, which is possible
becuase the integral operator is just a Fourier transform in this case.

To prove Theorem 1.1, we would like to apply the criterion of Theorem 1.8.
This leads to the consideration of the operator $T_{2}$
given by
\begin{equation}
(T_{2}f)(\lambda) = \frac{i\chi (a,b) (\lambda)}{2\sqrt{\lambda}} \int\limits_{0}^{\infty}
\exp \left(-2i\sqrt{\lambda} x +\frac{i}{\sqrt{\lambda}}\int\limits_{0}^{x}V(t)\,dt\right)f(x)\,dx,
\end{equation}
where $(a,b) \subset (0, \infty)$ and $a > 0,$ $b < \infty.$
We seek to show that the operator $T_{2}$ satisfies the $L_{2}-L_{2}$ bound
for every choice of $a,$ $b$ (although the value of the constant in the estimate may of course depend on this choice).
Theorem 1.1 will then follow immediately from Theorem 1.8 given that $(a,b)$ is an 
arbitrary proper subinterval of $(0, \infty).$

For the proof of the $L_{2}-L_{2}$ estimate  it suffices to assume  that 
the potential $V$ satisfies $|V(x)| \leq C(1+x)^{-\frac{1}{2}-\epsilon},$
with some $\epsilon >0.$
Let $2\sqrt{\lambda}=k.$ It is clear that it is sufficient to show the $L_{2}-L_{2}$ for the operator $T_{2}'$ given by
\begin{equation}
(T_{2}'f)(k)= \chi (a,b) \int\limits_{0}^{\infty}
\exp \left(-ik x +\frac{2i}{k}\int\limits_{0}^{x}V(t)\,dt\right)f(x)\,dx,
\end{equation}
The operator $T'_{2}$ looks like a pseudodifferential operator (restricted to the interval $(a,b)$) with a symbol
\[ a(k,x) =  \exp \left(\frac{i}{k}\int\limits_{0}^{x}
V(t)\,dt \right). \]
We remind that a symbol $a(k,x)$ belongs to an exotic class $S_{\rho, \sigma}$
if $a(k,x)$ is an infinitely differentiable function satisfying
\begin{equation}
 |\partial^{n}_{k}\partial^{m}_{x}a(k,x)| \leq C_{mn} (1+|x|)^{\sigma n
-\rho m}
\end{equation}
for every $m,$ $n.$
For the symbol classes $S_{\rho, \sigma},$ $1>\rho \geq \sigma \geq 0,$
the $L_{2}-L_{2}$ estimate is well-known (see, e.g., \cite{Stein1}). However,
for our purpose, although we may without loss of generality assume
that $V \in C^{\infty}$ (absorbing all lack of smoothness into short range
correction which is easy to treat), there is no hope in general that an estimate
like (29) holds for all integer $m,$ $n.$ Already taking the second derivative in $x,$ we should derivate $V,$ while under our assumptions we have absolutely no control over its derivative. However, Coifman and Meyer \cite{CoMe} have studied the question what is the minimal number of derivative estimates in (29) one has to ask for in order to have an $L_{2}-L_{2}$ bound. In particular, from their results it follows (Theorem
7 on page 30) that it suffices to check (29) for $m,$ $n=0,1$ for some
$1> \rho \geq \sigma \geq 0$ in order to ensure an $L_{2}-L_{2}$ bound
on $T.$ It is straightforward to check that for our symbol we have these
estimates for every $V$ satisfying $|V(x)| \leq C(1+x)^{-\frac{1}{2}-\epsilon}$
(and hence in particular for every $V$ satisfying  $|V(x)| \leq C(1+x)^{-\frac{2}{3}-\epsilon}):$
\[ |\partial_{x} a(k,x)| \leq C(1+x)^{-\frac{1}{2}
-\epsilon}; \]
\[ |\partial_{k} a(k,x)| \leq   \frac{2C}{(1-2\epsilon)a}
(1+x)^{\frac{1}{2}-\epsilon}; \]
\[| \partial_{x}\partial_{k} a(k,x)| \leq C_{11} (1+x)^{-2\epsilon}. \]
In particular, $a(k,x) \in ``S_{\frac{1}{2}, \frac{1}{2}}"$ (and even
$``S_{\frac{1}{2}+\epsilon, \frac{1}{2}-\epsilon}"$) where quotations mean the reduced
number of conditions on the derivatives, i.e. $m$ and $n$ are not
greater than $1$ in (29). Hence we have a theorem: \\

\noindent \bf Theorem 3.1. \it The operator $T_{2},$ given by \rm (27), \it satisfies
the $L_{2}-L_{2}$ bound $\|T_{2}f\|_{2} \leq C_{2}\|f\|_{2}$ if $V(x)$ verifies
$|V(x)| \leq C(1+x)^{-\frac{1}{2}-\epsilon}.$ Also, the bound $\|T_{2}f\|_{\infty} \leq C_{1}\|f\|_{1}$
holds trivially. \\

\rm This theorem together with criterion given by Theorem 1.8  implies
Theorem 1.1.

\rm We also sketch an alternative proof of Theorem 3.1 which
uses only an $L_{2}-L_{2}$ bound for the usual ``exotic" symbol class
$S_{\frac{1}{2}, \frac{1}{2}}$ with (29) true for any number of derivatives.
For this we need the following lemma: \\

\noindent \bf Lemma 3.2. \it Let $V(x)$ satisfy $|V(x)| \leq C(1+x)^{-\frac{1}{2}
-\epsilon};$ then we can represent a function $V(x)$ as a sum $V(x)=V_{1}(x)+
V_{2}(x),$ where $V_{1}(x)$ satisfies
\[ |V_{1}^{(m)}(x)|\leq C_{m}(1+x)^{-\frac{1}{2}(m+1)-\epsilon} \]
for every integer $m \geq 0,$ and $V_{2}(x)$ is conditionally integrable:
$\int_{0}^{x}V(t)\,dt$ converges as $x \rightarrow \infty.$ \\

\bf Proof. \rm Define an increasing sequence $\{a_{n} \}_{n=1}^{\infty}$ by
the conditions $a_{0} =1,$ $a_{n}-a_{n-1}=a_{n-1}^{\frac{1}{2}}.$ Let
$\xi$ be a $C^{\infty}$ function such that $\xi$ vanishes on the
interval $(-\delta, \delta),$ $\delta$ small positive number, and $\xi=1$
outside $(-2\delta, 2\delta).$ Let us define $V_{1}(x)$ by
\[ V_{1}(x) = \sum\limits_{n=1}^{\infty}C_{n} \chi (a_{n}, a_{n+1})
\xi \left(\frac{x-a_{n}}{a_{n}^{\frac{1}{2}}}\right)\xi \left(
\frac{x-a_{n+1}}{a_{n+1}^{\frac{1}{2}}} \right). \]
We choose each $C_{n}$ by a condition that $\int_{a_{n}}^{a_{n+1}}
(V-V_{1})(t)\,dt=0$ for every $n.$ It is easy to check that $V_{1}
\in C^{\infty}$ and for $x \in (a_{n}, a_{n+1})$ we have
\[ |V_{1}^{(m)}(x)| \leq \sup_{x\in(a_{n},a_{n+1})} |V(x)| C_{\xi}
a_{n}^{-\frac{1}{2}m},\] where $C_{\xi}$ depends only on $L_{\infty}$
norms of the derivatives of $\xi$ up to the $m$-th order. It is
easy to see that $\frac{a_{n+1}}{a_{n}} \stackrel{n
\rightarrow \infty}{\longrightarrow}0$ and hence we obtain $|V_{1}^{(m)}(x)| \leq
C_{m} (1+x)^{-\frac{1}{2}(m+1) -\epsilon}.$ On the other hand,
$\int^{a_{n}}_{0} V_{2}(t)\,dt =0$ for every $n$ and therefore it is
easy to see that $V_{2}$ is conditionally integrable and in fact
$|\int^{\infty}_{x} V_{2}(t)\,dt|\leq Cx^{-\epsilon}.$ $\Box$

Now let us write the symbol $a(k,x)$ as follows:
\[ a(k,x) = i\chi (a,b) (k)\exp \left(\frac{i}{2k}\int\limits_{0}^{x}V_{1}(t)
\,dt \right) \exp \left( \frac{i}{k} \int\limits_{0}^{x}V_{2}(t) \,dt \right). \]
The first two factors constitute a symbol from the $S_{\frac{1}{2}, \frac{1}{2}}$
class by inspection. Denote this symbol by $a_{1}(k,x).$ For every bounded
function $f$ of compact support a pseudodifferential operator $T_{1},$
associated  with the symbol $a_{1}(k,x),$ satisfies
$\|T_{1}f\|_{2} \leq C_{1}\|f\|_{2}.$  Since $V_{2}$ is a
conditionally integrable function,  there exists a constant $C_{2}$
such that $|\int_{0}^{x}V_{2}(t)\,dt| \leq C_{2}$ for every $x.$
Write the action of the operator $T$ as
\[ Tf(k)= T_{1} \left( \sum\limits^{\infty}_{j=1} \frac{1}{j!}
\left( \frac{i}{k} \int\limits_{0}^{x} V_{2}(t)\,dt \right) ^{j}
f(x) \right) = \]
\[ = \sum\limits_{j=1}^{\infty} \frac{1}{j!}\left(\frac{i}{k}\right)^{j}
T_{1} \left( \left( \int\limits_{0}^{x}V(t)\,dt\right)^{j}f(x)\right), \]
where the change of the orders of the action of $T_{1}$ and summation
is justified by the absolute convergence of the series. Hence, for 
every bounded function of compact support, we have
\[ \| Tf\|_{2} \leq C_{1}\sum\limits_{j=1}^{\infty}\frac{1}{j!a^{j}}C_{2}^{j}
\|f\|_{2} \leq C_{1} \exp (\frac{C_{2}}{a}) \|f\|_{2}. \,\,\,\, \Box \]

We now consider slowly decaying perturbations of Schr\"odinger
operators with periodic potentials. Let $U(x)$ be a periodic, piecewise
continuous function of period $T.$ It is a well-known fact
(see, e.g., \cite{ReSi}) that the spectrum of the operator
\[ \tilde{H}_{U}= -\frac{d^{2}}{dx^{2}}+U(x) \]
acting on $L^{2}(-\infty, \infty)$ is purely absolutely continuous of
multiplicity two and consisits of bands $[a_{n}, b_{n}],$ n=1,..., where
$a_{n}<b_{n} \leq a_{n+1}$ for every $n.$ First we will consider slowly
decaying perturbations for the case of an operator $H_{U},$ $U(x+T)=U(x)$
for every $x>0,$ defined on the semi-axis with some boundary condition at zero.
It is easy to see that the absolutely continuous spectrum of this operator
is of multiplicity one and coincides as a set with the absolutely continuous
spectrum of the corresponding whole-axis operator. This follows, for example,
from the existence of the Bloch solutions, which we will discuss shortly.

To prove the stabilty of the absolutely continuous spectrum of periodic
Schr\"odinger operators under a new class of slowly decaying perturbations, 
we would like to apply Theorem 1.8. For this we need to establish an $L_{2}-L_{2}$ bound
for an appropriate operator $T_{2}.$  Rather detailed
knowledge of the properties of the solutions $\theta (x, \lambda),$ which we 
choose to be the Bloch functions, is 
important to achieve this goal.

Let us recall the basic facts about the spectrum and the eigenfunctions of the
one-dimensional Schr\"odinger operators with periodic potentials; for the
missing proofs we refer to \cite{ReSi}.

For every band $[a_{n}, b_{n}]$ there exists a real analytic function
$\gamma (\lambda),$ which is called quasimomentum, such that $\gamma (\lambda)$
changes monotonically on $(a_{n}, b_{n})$ from $0$ to $\pi$ if $n$ is odd
, and from $\pi$ to $0$ if $n$ is even. The derivative $\gamma'(\lambda)$
might only vanish at the points $a_{n}$ or $b_{n}$ and in this case
respectively $b_{n-1}=a_{n}$ or $b_{n}=a_{n+1},$ i.e. there is no gap
between the bands. For every energy $\lambda \in (a_{n}, b_{n})$ there exists a
solution $\theta (x, \lambda),$ which is called a Bloch function, such that
\[ \theta (x+T, \lambda) = \exp (i\gamma(\lambda)) \theta (x, \lambda) \]
and \[ \theta '(x+T, \lambda) = \exp (i\gamma(\lambda))\theta'(x, \lambda). \]
The following lemma holds: \\
\noindent \bf Lemma 3.3. \it For every $\lambda \in (a_{n}, b_{n}),$ the
solutions $\theta (x, \lambda)$ and $\overline{\theta}(x, \lambda)$ are
linearly independent. \\

\noindent \bf Proof. \rm Indeed, suppose that $\overline{\theta}(x, \lambda)=
c\theta(x, \lambda);$ then we must have
\[ \overline{\theta}(x+T, \lambda)=
\exp (-i\gamma(\lambda))\overline{\theta}(x, \lambda)\]
 and \[ \overline{\theta}
(x+T, \lambda)= c\theta(x+T, \lambda)=c \exp (i\gamma(\lambda))\theta (x, \lambda). \]
Together this implies $\sin \gamma (\lambda)=0,$ which is not possible when
$\lambda \in (a_{n}, b_{n}).$ $\Box$\\

\noindent Hence, the Wronskian  $W[\theta, \overline{\theta}]= \Im(\theta\overline{\theta}')\neq 0$
when $\lambda \in (a_{n}, b_{n}).$ Next, we remind that the function $\theta
(x, \lambda),$ normalized by a condition $\|\theta(x, \lambda) \|_{L_{2}(0,T)}
=1,$ is real analytic in $\lambda$ as a function in $L_{2}(0,T)$ when $\lambda$
belongs to $(a_{n}, b_{n}).$
Moreover, we have \\
\bf Lemma 3.4. \it The solution $\theta (x, \lambda)$ is analytic in
$\lambda$ when $\lambda \in (a_{n}, b_{n})$ for every fixed $x \in [0,T].$
Moreover, the functions $\theta (x, \lambda),$ $\partial_{x} \theta
(x, \lambda),$ $\partial_{\lambda}\theta (x, \lambda)$ and $\partial^{2}_{x \lambda}\theta (x, \lambda)$
are continuous functions in every rectangle $[a'_{n}, b'_{n}]
 \times [0,T],$ where $a_{n}<a'_{n}<b'_{n}<b_{n}.$ \\
\noindent \bf Proof. \rm Let us consider two solutions of (4),
$y_{1}(x, \lambda)$ and $y_{2}(x, \lambda),$ satisfying $y_{1}(x, \lambda)=0,$
$y_{1}'(x, \lambda) =1$ and $y_{2}(x, \lambda)=1,$ $y_{2}'(x, \lambda)=0.$
Functions $y_{1},$ $y_{2}$ satisfy the properties claimed for $\theta$ in the
lemma by standart calculations using integral equations. Let
us represent the function $\theta (x, \lambda)$ as a linear combination of
these functions:
\[ \theta (x, \lambda) = c_{1}(\lambda)y_{1}(x, \lambda)+c_{2}(\lambda)
y_{2}(x, \lambda). \]
Consider now a vector
\[ v(\lambda)= y_{2}(x, \lambda) - \frac{\langle y_{1}(x, \lambda), y_{2}(x,
\lambda) \rangle_{L_{2}(0,T)}}{\|y_{1}(x, \lambda)\|_{L_{2}(0,T)}}y_{1}(x, \lambda) \]
for $\lambda \in [a_{n},b_{n}].$ This is an analyitc in $L_{2}(0,T)$ vector
with norm bounded away from zero on $[a_{n}, b_{n}]$ (we remind that solutions
$y_{1}(x, \lambda)$ and $y_{2}(x, \lambda)$ satisfy different boundary
conditions and their derivatives are in $x$ are bounded by some constant
in $[a_{n}, b_{n}] \times [0,T]$). We have
\[ \langle \theta (x, \lambda), v(\lambda) \rangle_{L_{2}(0,T)} = c_{2}(\lambda)
\left( \|y_{2} \|_{L_{2}(0,T)} - \frac{|\langle y_{1}, y_{2} \rangle_{L_{2}(0,T)}|^{2}}
{\|y_{1}\|_{L_{2}(0,T)}} \right). \]
Hence
\[ c_{2}(\lambda) = \frac{\langle \theta (x, \lambda), v(\lambda) \rangle_{L_{2}(0,T)}}
{\|v(\lambda)\|_{L_{2}(0,T)}^{2}}. \]
The denominator of the last expression is bounded away from zero and all
functions on the right hand side are real analytic when $\lambda \in (a_{n}, b_{n}),$
hence, $c_{2}(\lambda)$ is analytic in this interval. Similarly, we show the
analyticity of $c_{1}(\lambda).$ The statement of the lemma now follows
from the properties of $c_{1}(\lambda),$
$c_{2}(\lambda),$ $y_{1}(x, \lambda)$ and $y_{2}(x, \lambda).$ $\Box$ \\

\noindent \bf Proof of Theorem 1.2. \rm Let us consider some band $[a_{n}, b_{n}].$ Pick
an arbitrary interval $(a'_{n}, b'_{n}) \subset [a_{n}, b_{n}],$ such that
$a_{n} <a'_{n},$ $b'_{n} < b_{n}.$ This interval will serve us as a set $S_{i}$
from Theorem 1.8. To prove the Theorem 1.2 it suffices to show that an operator
$T_{2}$ defined by
\[ (T_{2}f)(\lambda) = \frac{\chi(a'_{n}, b'_{n})}{2\Im (\theta \overline{\theta}')}
\int\limits_{0}^{\infty} \overline{\theta}^{2}(x, \lambda) \exp\left(-\frac
{2}{\Im(\theta \overline{\theta}')}\int\limits_{0}^{x}V(t)|\theta(t,\lambda)|^{2}
\,dt \right) f(x)\,dx \]
satisfies the bound $\| T_{2}f\|_{L_{2}(a'_{n}, b'_{n})} \leq C\|f\|_{2}$ for
every bounded function $f$ of compact support. 

First we will show that the $L_{2}-L_{2}$ bound holds for an operator $T_{1}$
defined by (13):
\[ (T_{1}f)(\lambda)= \frac{\chi(S)}{\Im (\theta \overline{\theta}')} \int\limits_{0}^{\infty}
\overline{\theta}(x,\lambda)^{2} f(x)\,dx. \]
This will only prove that the absolutely continuous spectrum of Schr\"odinger
operators with periodic potentials is stable under perturbations $V(x)$
satisfying $|V(x)|\leq C(1+x)^{-\frac{3}{4}-\epsilon}.$ However, later we will see
that it is easy to adapt the proof to obtain the $L_{2}-L_{2}$ bound
for the operator $T_{2}.$

From the discussion of the properties of Bloch functions it follows that
the Wronskian $W[\theta, \overline{\theta}]=
\Im (\theta \overline{\theta}')$ is bounded away from zero
on $[a'_{n}, b'_{n}].$ Indeed, the Wronskian is continuous (and, in fact, 
real analytic) inside each band and vanishes only at $a_{n}$ or $b_{n}$ by
Lemma 3.3. Let us denote
\[ \omega_{n} =\inf_{\lambda \in (a'_{n}, b'_{n})} \Im (\theta \overline{\theta}'). \]
Also, the module of the derivative of the quasimomentum, $|\gamma'(\lambda)|,$
is bounded away from zero on $(a'_{n}, b'_{n}).$ Let
\[ \eta_{n} = \inf_{\lambda \in (a'_{n}, b'_{n})} |\theta'(\lambda)|. \]
Next, we note that the function
\[ \sigma (x, \lambda) = \left( \exp (-i\gamma(\lambda)
\frac{x}{T}) \theta (x, \lambda) \right)^{2} \]
is a periodic function with period $T.$ Let us consider the Fourier series for
$\sigma (x, \lambda):$
\[ \sigma (x, \lambda) = \sum\limits_{j} \exp \left( 2\pi ij\frac{x}{T}\right)
 \hat{\sigma}_{j}
(\lambda). \]
By $\hat{f}_{j}$ or $\hat{f}(k)$ we denote the Fourier transform of the
function $f$ in the discrete and continuous case respectively.
From the properties of $\theta (x, \lambda),$ it follows that $\partial_{x}
\sigma (x, \lambda)$ is a continuous function on $[a'_{n}, b'_{n}] \times [
0,T];$ let us denote
\[ \sigma_{n} = \sup_{[a'_{n}, b'_{n}] \times [0,T]} |\partial_{x} \sigma
(x, \lambda)|. \]
Now note that
\[ \|T_{1}f\|_{L_{2}(a'_{n}, b'_{n})} \leq \frac{1}{2\omega_{n}}\|\chi(a'_{n}, b'_{n})
\int\limits_{0}^{\infty} \exp (2i \gamma (\lambda) \frac{x}{T}) \sigma (x, \lambda)
f(x)\,dx \|_{L_{2}(a'_{n}, b'_{n})} = \]
\[ = \frac{1}{2\omega_{n}}\|
\int\limits_{0}^{\infty} \exp (2i \gamma (\lambda) \frac{x}{T}) \sum\limits_{j}
\left( \exp (2\pi ij\frac{x}{T} \hat{\sigma}_{j}(\lambda\right)
f(x)\,dx \|_{L_{2}(a'_{n}, b'_{n})}. \]
Since $|\partial\sigma (x, \lambda)| \leq \sigma_{n}$ for all $x \in
[0,T]$ uniformly in $\lambda \in (a'_{n}, b'_{n}),$ it is a standart fact
that the Fourier series (in $x$) for $\sigma (x, \lambda)$ converges absolutely.
In fact, even for Lipshitz-continuous with power $\alpha > \frac{1}{2}$
function $f(x)$ one has $\sum_{n} |\hat{f}(n)| \leq C\|f\|_{\Lambda_{\alpha}},$
see, for example, \cite{Kat}.

Hence, we can change the order of summation and integration in the previous
formula. We have
\[ \|T_{1}f\|_{L_{2}(a'_{n}, b'_{n})} \leq \frac{1}{2\omega_{n}}\|
\sum\limits_{j} \hat{\sigma}_{j}(\lambda)\int\limits_{0}^{\infty} \exp \left( 2i
\frac{x}{T}) (\gamma (\lambda) + j\pi) \right)
f(x)\,dx \|_{L_{2}(a'_{n}, b'_{n})} \leq \]
\[ \frac{1}{2\omega_{n}}\left(\, \int\limits_{a'_{n}}^{b'_{n}} \left( \sum_{j} |\hat{f} \left(
\frac{2(\gamma(\lambda)+j\pi)}{T} \right) \hat{\sigma}_{j}(\lambda) | \right)^{2}
\, d\lambda \right)^{\frac{1}{2}}  \leq \]
\[ \leq \frac{1}{2\omega_{n}} \left( \int\limits_{a'_{n}}^{b'{n}}
\left( \sum\limits_{j} |\hat{\sigma}_{j}|^{2}(\lambda)\right)\left( \sum\limits
_{j} \left|\hat{f}\left( \frac{2(\gamma(\lambda) +j\pi)}{T}\right)\right|^{2} \right)
\,d\lambda \right)^{\frac{1}{2}} \leq \]
\[ \leq \frac{\|\sigma^{2}(x,\lambda)\|_{L_{2}(0,T)}T^{\frac{1}{2}}}{2\omega_{n}\eta_{n}^{\frac{1}{2}}}
\left(\sum\limits_{j} \left|\,\int\limits_{j\pi +
\gamma (a'_{n})}^{j\pi +\gamma (b'_{n})} |\hat{f}(y)|^{2} \,dy\right|
\right)^{\frac{1}{2}}. \]
To obtain the last inequlity  we changed the orders
of summation and integration and introduced for each $j$ a new variable
$y = \frac{2(\gamma(\lambda)+j\pi)}{T}.$ We also note that from Lemma 3.4 it follows that 
$\sup_{\lambda \in [a_{n}',b_{n}']}\|\sigma^{2}(x, \lambda)\|_{L_{2}(0,T)} \leq
C_{n} <\infty.$ 
Hence, the last expression we obtained is estimated by
\[ \frac{C_{n}T^{\frac{1}{2}}}{2\omega_{n} \eta_{n}^{\frac{1}{2}}} \|\hat{f}\|_{2} \]
since the function $\gamma$ maps the interval $(a'_{n}, b'_{n})$ into
the interval $(0, \pi).$

Therefore, we get the desired bound
\[ \|T_{1}f\|_{L_{2}(a'_{n}, b'_{n})} \leq C \|f\|_{2}. \]
Now note that we can write the action of $T_{2}$ in a way similar to that
of $T_{1}:$
\begin{equation}
 (T_{2}f)(\lambda) = \frac{\chi (a_{n}', b_{n}')}{W [\theta, \overline{\theta}]}
\sum\limits_{j} \hat{\sigma}_{j}(\lambda) \int\limits_{0}^{\infty} \exp
\left(2i \frac{x}{T}(\gamma (\lambda) +j\pi)\right) a(\lambda, x) f(x)\,dx,
\end{equation}
where
\[ a(\lambda, x) = \exp \left( \frac{1}{W[\theta, \overline{\theta}]}\int\limits_{0}^{x} |\sigma(t, \lambda)|^{2}V(t)\,dt \right). \]
Proceeding with the estimation of the $L_{2}$ norm of the right-hand side of
(30) exactly as we did it before, we find that to establish the $L_{2}-L_{2}$
bound, it is sufficient to show that it holds for an operator $\tilde{T},$
defined by
\[ (\tilde{T}f)(y) = \int\limits_{0}^{\infty} \exp (iyx) \tilde{a}(y,x)\,dx, \]
with a ``symbol" $\tilde{a}(y,x)$ defined by 
\[ \tilde{a}(y,x)= \chi (2\gamma (a_{n}'), 2\gamma (b_{n}')) \times \] \[ \times \exp\left(\frac{1}{W[
\theta  (x,\gamma^{-1}(yT/2)), \overline{\theta}(x,\gamma^{-1}(yT/2)]}
\int\limits_{0}^{x}|\theta(t, \gamma^{-1}(yT/2))|^{2}V(t)\,dt\right) \] 
if $y \in [0, \frac{2\pi}{T}],$ 
and periodic in $y:$ $\tilde{a}(y+\frac{2\pi}{T},x) = \tilde{a}(y,x).$
The operator $\tilde{T}$ replaces the Fourier transform which appeared in the
estimate of $T_{1}.$ We note that the discontinuity in $y$ due to the presence of characteristic functions is artificail. We can  always replace the characteristic functions by smooth functions of compact support equal to $1$ when $y \in (2\gamma(a_{n}', \gamma(b_{n}'))$ and vanishing outside $(0, \pi).$ Form the $L_{2}-L_{2}$ bound for such operator would follow the bound for the original one. Now it is straightforward to check, using Lemma 3.4 and properties
of the quasimomentum $\gamma (\lambda),$ that we have
\[ |\partial_{x}^{\alpha}\partial_{y}^{\beta} \tilde{a}(y, x)| \leq
C_{\alpha \beta}(1+x)^{(-\frac{1}{2}-\epsilon)\alpha + (\frac{1}{2}-\epsilon)\beta} \]
for all $\alpha,$ $\beta$ taking values in $\{0,1\}.$ Hence, the "symbol''
$\tilde{a}(\lambda, x)$ belongs to the "$S_{\frac{1}{2},\frac{1}{2}}$'' class with the reduced number of
conditions on derivatives. By the Coifman-Meyer criterion \cite{CoMe}  it follows
that the operator $\tilde{T}$ satisfies an $L_{2}-L_{2}$ bound and therefore
this bound also holds for $T_{2}.$ 
  $\Box$ \\

To prove the whole-axis analog of Theorem 1.2, Theorem 1.4, we apply the whole
axis criterion formulated in Theorem 1.10. The needed $L_{2}-L_{2}$ bounds
are obtained similarly to the semi-axis case.

As a final remark we note that the results parallel to those we show here
also hold for Jacobi matrices case. The role of key Theorem 2.1 is played
by its discrete analog ( which in particular follows from considerations
in \cite{Kis1}). We plan to further develop this  theme in 
a subsequent publication.

\begin{center} 
\section*{ Acknowledgment }
\end{center}

\normalsize

I would like to thank Prof. B.~Simon for stimulating discussions and valuable comments. I am very grateful to Prof. S.~Semmes for inspiring and 
informative conversations on harmonic analysis and to Prof. F.~Gesztesy
for asking me questions which largely motivated this work.

I gratefully acknowledge hospitality of IHES, where part of 
this work was done. Research at  MSRI supported in part by NSF grant DMS 
9022140.

\end{document}